\documentclass[12pt,amssymb,amstex]{article}

\usepackage{amsmath}
\usepackage{amssymb}
\usepackage{amsxtra}
\usepackage{latexsym}
\usepackage{ifthen}

\usepackage[dvips]{graphicx}
\usepackage{epstopdf}
\usepackage{epsfig}
\usepackage{epic}
\usepackage{eepic}

\def\Xint#1{\mathchoice
   {\XXint\displaystyle\textstyle{#1}}
   {\XXint\textstyle\scriptstyle{#1}}
   {\XXint\scriptstyle\scriptscriptstyle{#1}}
   {\XXint\scriptscriptstyle\scriptscriptstyle{#1}}
   \!\int}
\def\XXint#1#2#3{{\setbox0=\hbox{$#1{#2#3}{\int}$}
     \vcenter{\hbox{$#2#3$}}\kern-.5\wd0}}
\def\dashint{\Xint-}

\begin{document}

\newcounter{lemma}[section]
\newcommand{\lemma}{\par \refstepcounter{lemma}%
{\bf Lemma \arabic{section}.\arabic{lemma}.}}
\renewcommand{\thelemma}{\thesection.\arabic{lemma}}

\newcounter{corol}[section]
\newcommand{\corol}{\par \refstepcounter{corol}%
{\bf Corollary \arabic{section}.\arabic{corol}.}}
\renewcommand{\thecorol}{\thesection.\arabic{corol}}

\newcounter{rem}[section]
\newcommand{\rem}{\par \refstepcounter{rem}%
{\bf Remark \arabic{section}.\arabic{rem}.}}
\renewcommand{\therem}{\thesection.\arabic{rem}}

\newcounter{theo}[section]
\newcommand{\theo}{\par \refstepcounter{theo}%
{\bf Theorem \arabic{section}.\arabic{theo}.}}
\renewcommand{\thetheo}{\thesection.\arabic{theo}}

\newcounter{propo}[section]
\newcommand{\propo}{\par \refstepcounter{propo}%
{\bf Proposition \arabic{section}.\arabic{propo}.}}
\renewcommand{\thepropo}{\thesection.\arabic{propo}}

\numberwithin{equation}{section}

\newcommand{\osc}{\operatornamewithlimits{osc}}

\def\Xint#1{\mathchoice
   {\XXint\displaystyle\textstyle{#1}}%
   {\XXint\textstyle\scriptstyle{#1}}%
   {\XXint\scriptstyle\scriptscriptstyle{#1}}%
   {\XXint\scriptscriptstyle\scriptscriptstyle{#1}}%
   \!\int}
\def\XXint#1#2#3{{\setbox0=\hbox{$#1{#2#3}{\int}$}
     \vcenter{\hbox{$#2#3$}}\kern-.5\wd0}}
\def\dashint{\Xint-}

\def\cc{\setcounter{equation}{0}
\setcounter{figure}{0}\setcounter{table}{0}}

\overfullrule=0pt

\title{{\bf On boundary behavior and\\
Dirichlet problem for Beltrami equations}}

\author{\bf D. Kovtonyuk, I. Petkov,  V. Ryazanov and R. Salimov}

\date{\today \hskip 4mm ({\tt ARXIV-260112-KPRS.tex})}

\maketitle

\abstract We show that homeomorphic $W^{1,1}_{\rm loc}$ solutions
to the Beltrami equations $\overline{\partial}f=\mu\,\partial f$
satisfy certain moduli inequalities. On this basis, we develope
the theory of the boundary behavior of such solutions and prove a
series of criteria for the existence of regular, pseudoregular and
multi-valued solutions for the Dirichlet problem to the Beltrami
equations in Jordan domains and finitely  connected domains.
\endabstract

\medskip



\large

\section{Introduction}

Let $D$ be a domain in the complex plane ${\Bbb C}$, i.e., a
connected open subset of ${\Bbb C}$, and let $\mu:D\to{\Bbb C}$ be
a measurable function with $|\mu(z)|<1$ a.e. (almost everywhere) in
$D$. A {\bf Beltrami equation} is an equation of the form
\begin{equation}\label{eqBeltrami} f_{\bar z}=\mu(z)\,f_z\end{equation} where
$f_{\bar z}=\overline{\partial}f=(f_x+if_y)/2$,
$f_{z}=\partial f=(f_x-if_y)/2$, $z=x+iy$, and $f_x$ and $f_y$ are
partial derivatives of $f$ in $x$ and $y$, correspondingly. The
function $\mu$ is called the {\bf complex coefficient} and
\begin{equation}\label{eqKPRS1.1}K_{\mu}(z)=\frac{1+|\mu(z)|}{1-|\mu(z)|}\end{equation}
the {\bf dilatation quotient} of the equation (\ref{eqBeltrami}).
The Beltrami equation (\ref{eqBeltrami}) is said to be {\bf
degenerate} if ${\rm ess}\,{\rm sup}\,K_{\mu}(z)=\infty$. The
existence of homeomorphic $W^{1,1}_{\rm loc}$ solutions was
recently established to many degenerate Beltrami equations, see,
e.g., related references in the recent monograph \cite{MRSY} and
in the surveys \cite{GRSY} and \cite{SY}.

Boundary value problems for the Beltrami equations are due to the
well-known Riemann dissertation in the case of $\mu(z)=0$ and to the
papers of Hilbert (1904, 1924) and Poincare (1910) for the
corresponding Cauchy--Riemann system. The Dirichlet problem was well
studied for uniformly elliptic systems, see, e.g., \cite{Boj} and
\cite{Vekua}. The Dirichlet problem for degenerate Beltrami
equations in the unit disk was studied in \cite{Dy}. However, the
criteria for the existence of solutions for the Dirichlet problem in
\cite{Dy} are not invariant under conformal mappings of Riemann.
Hence we give here the corresponding theorems on the existence of
regular solutions in arbitrary Jordan domains as well as of
pseudoregular and multi-valued solutions in arbitrary multiply
connected domains  bounded by a finite collection of mutually
disjoint Jordan curves. In comparison with the work \cite{Dy}, our
approach is based on estimates of the modulus of dashed lines but
not of paths under arbitrary homeomorphic $W^{1,1}_{\rm loc}$
solutions of the Beltrami equations.

Recall that every holomorphic (analytic) function $f$ in a domain
$D\subset{\Bbb C}$ satisfies the simplest Beltrami equation
\begin{equation}\label{eqKPRS1.2}f_{\bar z}=0\end{equation} with
$\mu(z)\equiv0$. If a holomorphic function $f$ given in the unit
disk ${\Bbb D}$ is continuous in its closure, then by the Schwarz
formula \begin{equation}\label{eqKPRS1.3}f(z)=i\,{\rm
Im}\,f(0)+\frac{1}{2\pi i}\int\limits_{|\zeta|=1}{\rm
Re}\,f(\zeta)\cdot\frac{\zeta+z}{\zeta-z}\frac{d\zeta}{\zeta}\,,\end{equation}
see, e.g., Section 8, Chapter III, Part 3 in \cite{HuCo}. Thus,
the holomorphic function $f$ in the unit disk ${\Bbb D}$ is
determinated, up to a purely imaginary additive constant $ic$,
$c={\rm Im}\,f(0)$, by its real part $\varphi(\zeta)={\rm
Re}\,f(\zeta)$ on the boundary of ${\Bbb D}$.

\medskip

Hence the {\bf Dirichlet problem} for the Beltrami equation
(\ref{eqBeltrami}) in a domain $D\subset{\Bbb C}$ is the problem on
the existence of a continuous function $f:D\to{\Bbb C}$ having
partial derivatives of the first order a.e., satisfying
(\ref{eqBeltrami}) a.e. and such that
\begin{equation}\label{eqGrUsl}\lim\limits_{z\to\zeta}{\rm
Re}\,f(z)=\varphi(\zeta)\qquad\forall\ \zeta\in\partial
D\end{equation} for a prescribed continuous function
$\varphi:\partial D\to{\Bbb R}$. It is obvious that if $f$ is a
solution of this  problem, then the function $F(z)=f(z)+ic$,
$c\in{\Bbb R}$, is so.

\medskip

Here we show that the Dirichlet problem (\ref{eqGrUsl}) has regular
solutions in an arbitrary Jordan domain and pseudoregular and
multi-valued solutions in an arbitrary finitely connected domain for
wide classes of degenerate Beltrami equations (\ref{eqBeltrami}).

\medskip

Finally, it is necessary to note that the existence of solutions
for the Dirichlet problem to the degenerate Beltrami equations
with two characteristics \begin{equation}\label{eqBeltrami2}
f_{\bar z}=\mu(z)\,f_z+\nu(z)\,\overline{f_z}\end{equation}
remains open although the corresponding theorems on the existence
of homeomorphic $W_{\rm loc}^{1,1}$ solutions of
(\ref{eqBeltrami2}) have been established in the series of the
recent papers \cite{BGR$_1$}--\cite{BGR$_3$}. This problem is
important because the Beltrami equations of the second type
\begin{equation}\label{eqBeltrami3} f_{\bar z}=\nu(z)\,\overline{f_z}\end{equation}
play a great role in many problems of mathematical physics, see,
e.g., \cite{KK}. However, the Dirichlet problem for the equation
(\ref{eqBeltrami3}) demands an essential modification of our
approach.

\section{Preliminaries}

Throughout this paper, $B(z_0,r)=\{z\in{\Bbb C}:|z_0-z|<r\}$,
${\Bbb D}=B(0,1)$, $S(z_0,r)=\{z\in{\Bbb C}:|z_0-z|=r\}$,
$S(r)=S(0,r)$, $R(z_0,r_1,r_2)=\{z\in{\Bbb C}:r_1<|z-z_0|<r_2\}$.

The class BMO was introduced by John and Nirenberg (1961) in the
paper \cite{JN} and soon became an important concept in harmonic
analysis, partial differential equations and related areas; see,
e.g., \cite{HKM} and \cite{RR}.

Recall that a real-valued function $u$ in a domain $D$ in ${\Bbb
C}$ is said to be of {\bf bounded mean oscillation} in $D$, abbr.
$u\in{\rm BMO}(D)$, if $u\in L_{\rm loc}^1(D)$ and
\begin{equation}\label{lasibm_2.2}\Vert u\Vert_{*}:=
\sup\limits_{B}{\frac{1}{|B|}}\int\limits_{B}|u(z)-u_{B}|\,dm(z)<\infty\,,\end{equation}
where the supremum is taken over all discs $B$ in $D$, $dm(z)$
corresponds to the Lebesgue measure in ${\Bbb C}$ and
$$u_{B}={\frac{1}{|B|}}\int\limits_{B}u(z)\,dm(z)\,.$$ We write $u\in{\rm BMO}_{\rm loc}(D)$ if
$u\in{\rm BMO}(U)$ for every relatively compact subdomain $U$ of $D$
(we also write BMO or ${\rm BMO}_{\rm loc }$ if it is clear from the
context what $D$ is).

\medskip

Following the paper \cite{IR}, see also \cite{IRp} and
\cite{MRSY}, we say that a function $\varphi:D\to{\Bbb R}$ has
{\bf finite mean oscillation} at a point $z_0\in D$ if
\begin{equation}\label{FMO_eq2.4}\overline{\lim\limits_{\varepsilon\to0}}\ \ \
\dashint_{B(z_0,\varepsilon)}|{\varphi}(z)-\widetilde{\varphi}_{\varepsilon}(z_0)|\,dm(z)<\infty\,,\end{equation}
where \begin{equation}\label{FMO_eq2.5}
\widetilde{\varphi}_{\varepsilon}(z_0)=\dashint_{B(z_0,\varepsilon)}
{\varphi}(z)\,dm(z)\end{equation} is the mean value of the
function ${\varphi}(z)$ over the disk $B(z_0,\varepsilon)$. Note
that the condition (\ref{FMO_eq2.4}) includes the assumption that
$\varphi$ is integrable in some neighborhood of the point $z_0$.
We say also that a function $\varphi:D\to{\Bbb R}$ is of {\bf
finite mean oscillation in $D$}, abbr. $\varphi\in{\rm FMO}(D)$ or
simply $\varphi\in{\rm FMO}$, if $\varphi\in{\rm FMO}(z_0)$ for
all points $z_0\in D$. We write $\varphi\in{\rm
FMO}(\overline{D})$ if $\varphi$ is given in a domain $G$ in
$\Bbb{C}$ such that $\overline{D}\subset G$ and $\varphi\in{\rm
FMO}(z_0)$ for all $z_0\in\overline{D}$.

The following statement is obvious by the triangle inequality.

\medskip

\begin{propo}\label{FMO_pr2.1} {\it If, for a  collection of numbers
$\varphi_{\varepsilon}\in{\Bbb R}$,
$\varepsilon\in(0,\varepsilon_0]$,
\begin{equation}\label{FMO_eq2.7}\overline{\lim\limits_{\varepsilon\to0}}\ \ \
\dashint_{B(z_0,\varepsilon)}|\varphi(z)-\varphi_{\varepsilon}|\,dm(z)<\infty\,,\end{equation}
then $\varphi $ is of finite mean oscillation at $z_0$.}
\end{propo}

In particular choosing in Proposition \ref{FMO_pr2.1},
$\varphi_{\varepsilon}\equiv0$, $\varepsilon\in(0,\varepsilon_0]$,
we obtain the following.

\medskip

\begin{corol}\label{FMO_cor2.1} {\it If, for a point $z_0\in D$,
\begin{equation}\label{FMO_eq2.8}\overline{\lim\limits_{\varepsilon\to 0}}\ \ \
\dashint_{B(z_0,\varepsilon)}|\varphi(z)|\,dm(z)<\infty\,,
\end{equation} then $\varphi$ has finite mean oscillation at
$z_0$.} \end{corol}

\medskip

Recall that a point $z_0\in D$ is called a {\bf Lebesgue point} of
a function $\varphi:D\to{\Bbb R}$ if $\varphi$ is integrable in a
neighborhood of $z_0$ and \begin{equation}\label{FMO_eq2.7a}
\lim\limits_{\varepsilon\to 0}\ \ \ \dashint_{B(z_0,\varepsilon)}
|\varphi(z)-\varphi(z_0)|\,dm(z)=0\,.\end{equation} It is known
that, almost every point in $D$ is a Lebesgue point for every
function $\varphi\in L^1(D)$. Thus we have by Proposition
\ref{FMO_pr2.1} the following corollary.

\medskip

\begin{corol}\label{FMO_cor2.7b} {\it Every
locally integrable function $\varphi:D\to{\Bbb R}$, has a finite
mean oscillation at almost every point in $D$.} \end{corol}

\medskip

\begin{rem}\label{FMO_rmk2.13a} Note that the function $\varphi(z)=\log\left(1/|z|\right)$
belongs to BMO in the unit disk $\Delta$, see, e.g., \cite{RR}, p.
5, and hence also to FMO. However,
$\widetilde{\varphi}_{\varepsilon}(0)\to\infty$ as
$\varepsilon\to0$, showing that condition (\ref{FMO_eq2.8}) is
only sufficient but not necessary for a function $\varphi$ to be
of finite mean oscillation at $z_0$. Clearly, ${\rm
BMO}(D)\subset{\rm BMO}_{\rm loc}(D)\subset{\rm FMO}(D)$ and as
well-known ${\rm BMO}_{\rm loc}\subset L_{\rm loc}^p$ for all
$p\in[1,\infty)$, see, e.g., \cite{BN}, by FMO is not a subclass
of $L_{\rm loc}^p$ for any $p>1$ but only $L_{\rm loc}^1$. Thus
the class FMO is much more wide than ${\rm BMO}_{\rm
loc}$.\end{rem}

\medskip

Versions of the next lemma has been first proved for the class BMO
in the planar case in \cite{RSY$_1$}, \cite{RSY$_8$}, and then in
the space case in \cite{MRSY$_5$}, \cite{MRSY$_6$}. For the FMO
case, see the papers \cite{IR}, \cite{IRp}, \cite{RS},
\cite{RSY$_7$}, \cite{RSY$_4$},  and the monograph \cite{MRSY}.

\medskip

\begin{lemma}\label{lem13.4.2} {\it Let $D$ be a domain in ${\Bbb C}$ and let
$\varphi:D\to{\Bbb R}$ be a  non-negative function  of the class
${\rm FMO}(z_0)$ for some $z_0\in D$. Then
\begin{equation}\label{eq13.4.5}\int\limits_{\varepsilon<|z-z_0|<\varepsilon_0}\frac{\varphi(z)\,dm(z)}
{\left(|z-z_0|\log\frac{1}{|z-z_0|}\right)^2}=O\left(\log\log\frac{1}{\varepsilon}\right)\
\quad\text{as}\quad\varepsilon\to 0\end{equation} for some
$\varepsilon_0\in(0,\delta_0)$ where $\delta_0=\min(e^{-e},d_0)$,
$d_0=\sup_{z\in D}|z-z_0|$.} \end{lemma}

\medskip

Recall a connection between some integral conditions, see, e.g.,
\cite{RSY12}--\cite{RSY}.

\medskip

\begin{theo}\label{th5.555} {\it Let $Q:{\Bbb D}\to[0,\infty]$ be a
measurable function such that \begin{equation}\label{eq5.555}
\int\limits_{\Bbb D}\Phi(Q(z))\,dm(z)<\infty\end{equation} where
$\Phi:[0,\infty]\to[0,\infty]$ is a non-decreasing convex function
such that \begin{equation}\label{eq3.333a}
\int\limits_{\delta_0}^{\infty}\frac{d\tau}{\tau\Phi^{-1}(\tau)}=\infty\end{equation}
for some $\delta_0>\Phi(0)$. Then \begin{equation}\label{eq3.333A}
\int\limits_0^1\frac{dr}{rq(r)}=\infty\end{equation} where $q(r)$
is the average of the function $Q(z)$ over the circle $|z|=r$.}
\end{theo}

\medskip

The following lemma is also useful, see Lemma 2.1 in \cite{KR$_2$}
or Lemma 9.2 in \cite{MRSY}.

\medskip

\begin{lemma}\label{lem8.4.1} {\it Let $(X,\mu)$ be a measure space with a
finite measure $\mu$, $p\in(1,\infty)$ and let
$\varphi:X\to(0,\infty)$ be a measurable function. Set
\begin{equation}\label{eq8.4.2}I(\varphi,p)=\inf\limits_{\alpha}\int\limits_{X}\varphi\,\alpha^p\,d\mu\end{equation}
where the infimum is taken over all measurable functions
$\alpha:X\to[0,\infty]$ such that \begin{equation}\label{eq8.4.3}
\int\limits_{X}\alpha\,d\mu=1\,.\end{equation} Then
\begin{equation}\label{eq8.4.4}I(\varphi,p)=\left[\int\limits_{X}\varphi^{-\lambda}\,d\mu\right]^{-\frac{1}{\lambda}}\end{equation}
where \begin{equation}\label{eq8.4.5}\lambda=\frac{q}{p}\,,\qquad
\frac{1}{p}+\frac{1}{q}=1\,,\end{equation} i.e.
$\lambda=1/(p-1)\in(0,\infty)$. Moreover, the infimum in
(\ref{eq8.4.2}) is attained only for the function
\begin{equation}\label{eq8.4.6}\alpha_0=C\cdot\varphi^{-\lambda}\end{equation} where
\begin{equation}\label{eq8.4.7}C=\left(\int\limits_{X}\varphi^{-\lambda}\,d\mu\right)^{-1}\,.\end{equation}}
\end{lemma}

\medskip

Finally, recall that the {\bf (conformal) modulus} of a family
$\Gamma$ of paths $\gamma$ in ${\Bbb C}$ is the quantity
\begin{equation}\label{eqModul} M(\Gamma)=\inf_{\varrho\in{\rm adm}\,\Gamma}\int\limits_{\Bbb
C}\varrho^2(z)\,dm(z)\end{equation} where a Borel function
$\varrho:{\Bbb C}\to[0,\infty]$ is {\bf admissible} for $\Gamma$,
write $\varrho\in{\rm adm}\,\Gamma$, if
\begin{equation}\label{eqAdm}\int\limits_{\gamma}\varrho\,ds\geqslant1\quad\forall\ \gamma\in\Gamma\,.\end{equation}
Here $s$ is a natural parameter of the length on $\gamma$.

\section{On regular domains}

First of all, recall the following topological notion. A domain
$D\subset{\Bbb C}$ is said to be {\bf locally connected at a
point} $z_0\in\partial D$ if, for every neighborhood $U$ of the
point $z_0$, there is a neighborhood $V\subseteq U$ of $z_0$ such
that $V\cap D$ is connected. If this condition holds for all
$z_0\in \partial D$, then $D$ is said to be locally connected on
$\partial D$. Note that every Jordan domain $D$ in ${\Bbb C}$ is
locally connected on  $\partial D$, see, e.g., \cite{Wi}, p. 66.

\begin{figure}[h]
\centerline{\includegraphics[scale=0.25]{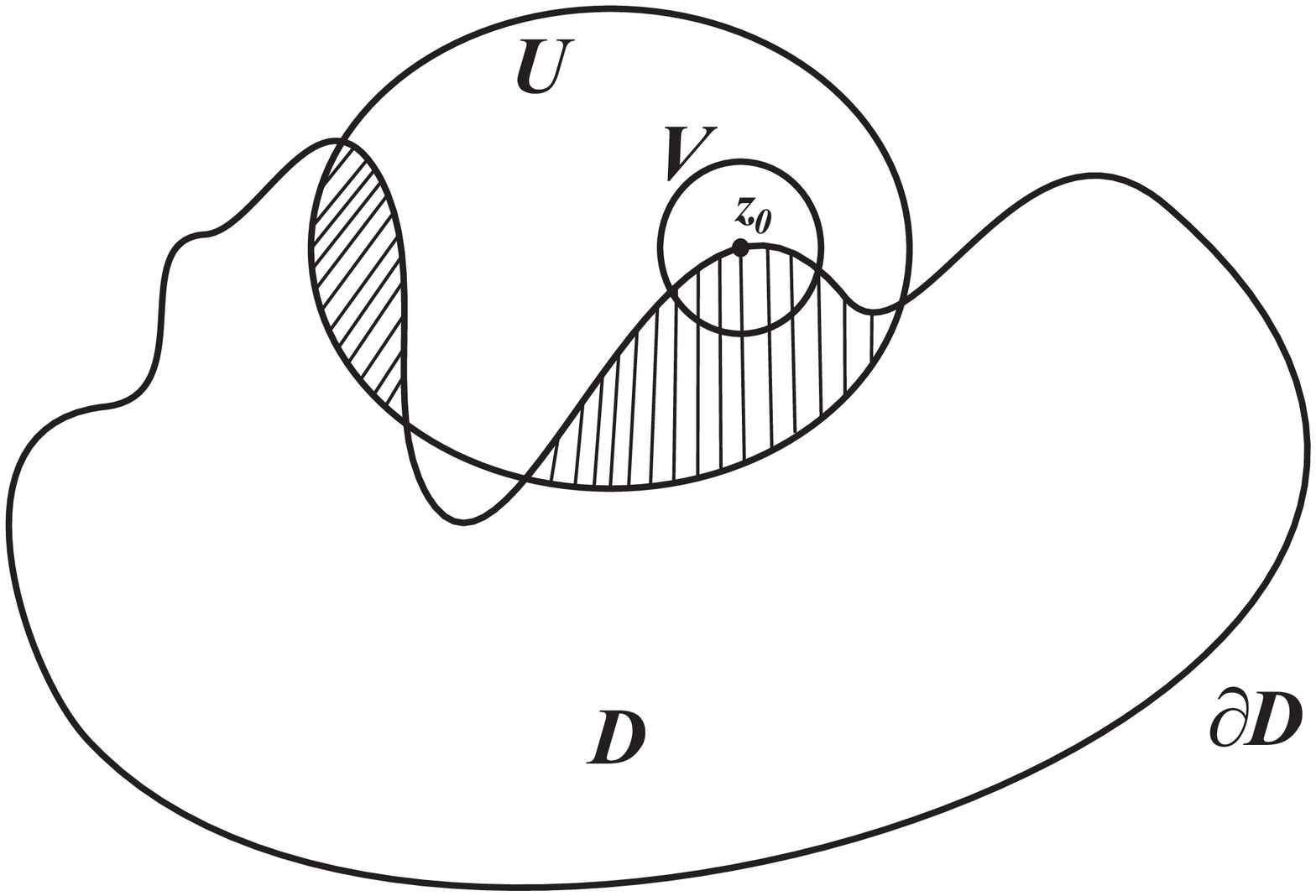}}
\end{figure}

We say that $\partial D$ is {\bf weakly flat at a point}
$z_0\in\partial D$ if, for every neighborhood $U$ of the point
$z_0$ and every number $P>0$, there is a neighborhood $V\subset U$
of $z_0$ such that \begin{equation}\label{eq1.5KR}
M(\Delta(E,F;D))\geqslant P\end{equation} for all continua $E$ and
$F$ in $D$ intersecting $\partial U$ and $\partial V$. We say that
$\partial D$ is {\bf weakly flat} if it is weakly flat at each
point $z_0\in\partial D$.

\begin{figure}[h]
\centerline{\includegraphics[scale=0.35]{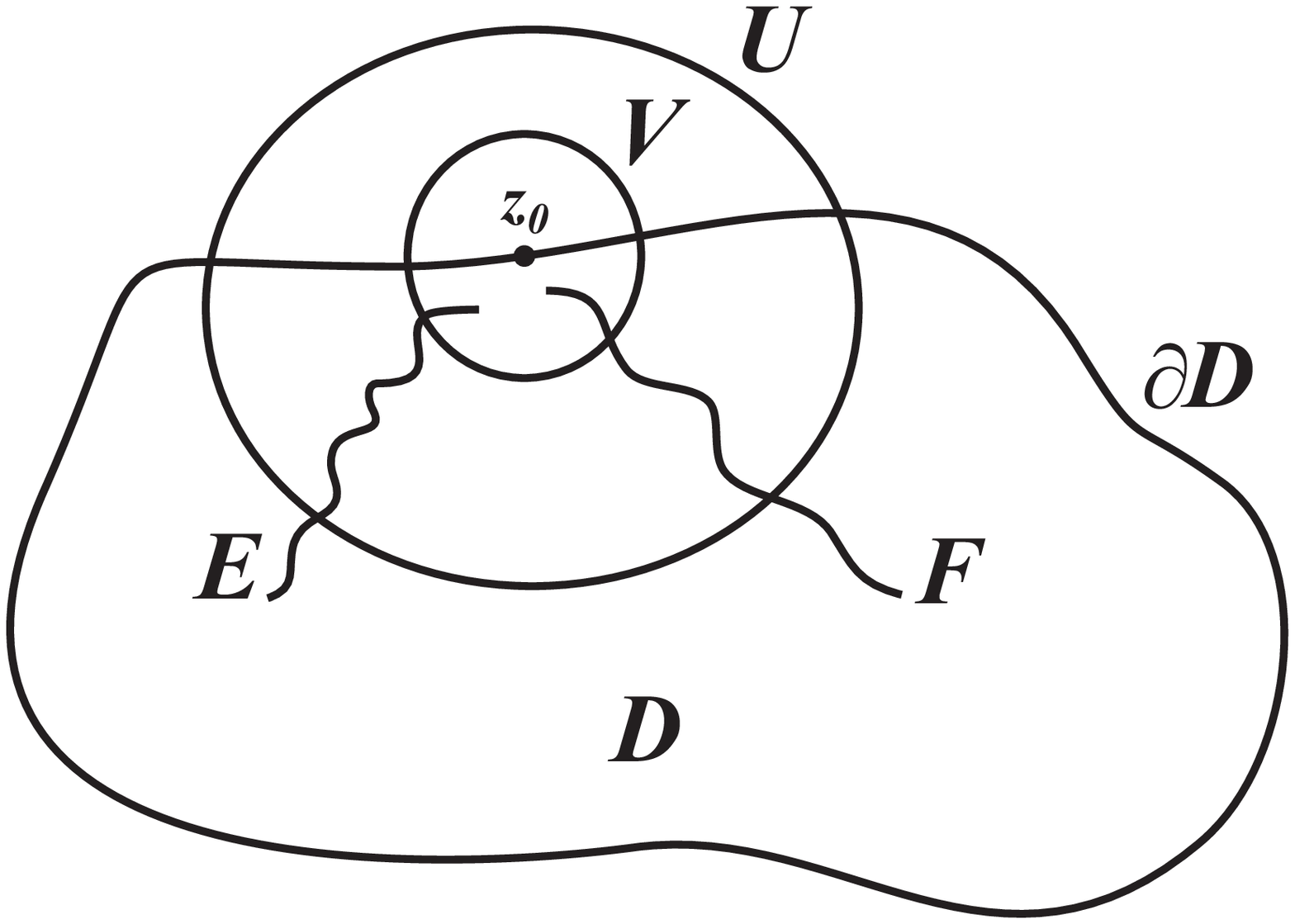}}
\end{figure}

We also say that a point $z_0\in\partial D$ is {\bf strongly
accessible} if, for every neighborhood $U$ of the point $z_0$,
there exist a compactum $E$ in $D$, a neighborhood $V\subset U$ of
$z_0$ and a number $\delta>0$ such that
\begin{equation}\label{eq1.6KR}M(\Delta(E,F;D))\geqslant\delta\end{equation}
for all continua $F$ in $D$ intersecting $\partial U$ and
$\partial V$. We say that $\partial D$ is {\bf strongly
accessible} if each point $z_0\in\partial D$ is strongly
accessible.

Here, in the definitions of strongly accessible and weakly flat
boundaries, one can take as neighborhoods $U$ and $V$ of a point
$z_0$ only balls (closed or open) centered at $z_0$ or only
neighborhoods of $z_0$ in another fundamental system of
neighborhoods of $z_0$. These conceptions can also be extended in a
natural way to the case of $\overline{\Bbb C}$ and $z_0=\infty$.
Then we must use the corresponding neighborhoods of $\infty$.

It is easy to see that if a domain $D$ in ${\Bbb C}$ is weakly flat
at a point $z_0\in\partial D$, then the point $z_0$ is strongly
accessible from $D$. Moreover, it was proved by us that if a domain
$D$ in ${\Bbb C}$ is weakly flat at a point $z_0\in\partial D$, then
$D$ is locally connected at $z_0$, see, e.g., Lemma 5.1 in
\cite{KR$_2$} or Lemma 3.15 in \cite{MRSY}.

The notions of strong accessibility and weak flatness at boundary
points of a domain in ${\Bbb C}$ defined in \cite{KR$_0$}, see also
\cite{KR$_2$} and \cite{RS}, are localizations and generalizations
of the corresponding notions introduced in \cite{MRSY$_5$} and
\cite{MRSY$_6$}, cf. with the properties $P_1$ and $P_2$ by
V\"ais\"al\"a in \cite{Va} and also with the quasiconformal
accessibility and the quasiconformal flatness by N\"akki in
\cite{Na$_1$}. Many theorems on a homeomorphic extension to the
boundary of quasiconformal mappings and their generalizations are
valid under the condition of weak flatness of boundaries. The
condition of strong accessibility plays a similar role for a
continuous extension of the mappings to the boundary.

A domain $D\subset{\Bbb C}$ is called a {\bf quasiextremal distance
domain}, abbr. {\bf QED-domain}, see \cite{GM}, if
\begin{equation}\label{e:7.1}M(\Delta(E,F;{\Bbb C})\leqslant K\cdot
M(\Delta(E,F;D))\end{equation} for some $K\geqslant1$ and all pairs
of nonintersecting continua $E$ and $F$ in $D$.

It is well known, see, e.g., Theorem 10.12 in \cite{Va}, that
\begin{equation}\label{eqKPR2.2}M(\Delta(E,F;{\Bbb C}))\geqslant\frac{2}{\pi}\log{\frac{R}{r}}\end{equation}
for any sets $E$ and $F$ in ${\Bbb C}$ intersecting all the
circles $S(z_0,\rho)$, $\rho\in(r,R)$. Hence a QED-domain has a
weakly flat boundary. One example in \cite{MRSY}, Section 3.8,
shows that the inverse conclusion is not true even in the case of
simply connected domains in ${\Bbb C}$.

A domain $D\subset{\Bbb C}$ is called a {\bf uniform domain} if
each pair of points $z_1$ and $z_2\in D$ can be joined with a
rectifiable curve $\gamma$ in $D$ such that
\begin{equation}\label{e:7.2} s(\gamma)\leqslant
a\cdot|z_1-z_2|\end{equation} and \begin{equation}\label{e:7.3}
\min\limits_{i=1,2}\ s(\gamma(z_i,z))\leqslant b\cdot{\rm
dist}(z,\partial D)\end{equation} for all $z\in\gamma$ where
$\gamma(z_i,z)$ is the portion of $\gamma$ bounded by $z_i$ and
$z$, see \cite{MaSa}. It is known that every uniform domain is a
QED-domain but there exist QED-domains that are not uniform, see
\cite{GM}. Bounded convex domains and bounded domains with smooth
boundaries are simple examples of uniform domains and,
consequently, QED-domains as well as domains with weakly flat
boundaries.

In the mapping theory and in the theory of differential equations,
it is also often applied the so-called Lipschitz boundaries.

They say that a domain $D$ in $\Bbb{C}$ is {\bf Lipschitz} if
every point $z_0\in\partial D$ has a neighborhood $U$ that can be
mapped by a bi-Lipschitz homeomorphism $\varphi$ onto the unit
disk ${\Bbb D}\subset{\Bbb C}$ in such a way that
$\varphi(\partial D\cap U)$ is the intersection of ${\Bbb D}$ with
the real axis. Note that a bi-Lipschitz homeomorphism is
quasiconformal and the modulus is a quasiinvariant under such
mappings. Hence the Lipschitz domains have weakly flat boundaries.

Recall that a map $\varphi:U\to{\Bbb C}$ is said to be {\bf
Lipcshitz} provided $|\varphi(z_1)-\varphi(z_2)|\leqslant
M\cdot|z_1-z_2|$ for some $M<\infty$ and for all $z_1$ and $z_2\in
U$, and {\bf bi-Lipcshitz} if in addition
$M^*|z_1-z_2|\leqslant|\varphi(z_1)-\varphi(z_2)|$ for some
$M^*>0$ and for all $z_1$ and $z_2\in U$.

\section{On estimates of modulus of dashed lines}

A continuous mapping $\gamma$ of an open subset $\Delta$ of the
real axis ${\Bbb R}$ or a circle into $D$ is called a {\bf dashed
line}, see, e.g., 6.3 in \cite{MRSY}. Note that such a set
$\Delta$ consists of a countable collection of mutually disjoint
intervals in ${\Bbb R}$. This is the motivation for the term. The
notion of the modulus of the family $\Gamma$ of dashed lines
$\gamma$ is defined similarly to (\ref{eqModul}). We say that a
property $P$ holds for {\bf a.e.} (almost every) $\gamma\in\Gamma$
if a subfamily of all lines in $\Gamma$ for which $P$ fails has
the modulus zero, cf. \cite{Fu}. Later on, we also say that a
Lebesgue measurable function $\varrho:{\Bbb C}\to[0,\infty]$ is
{\bf extensively admissible} for $\Gamma$, write $\varrho\in{\rm
ext\,adm}\,\Gamma$, if (\ref{eqAdm}) holds for a.e.
$\gamma\in\Gamma$, see, e.g., Section 9.2 in \cite{MRSY}.

\medskip

\begin{theo}\label{thDIR1} {\it Let $f$ be a homeomorphic $W^{1,1}_{\rm
loc}$ solution of the Beltrami equation (\ref{eqBeltrami}) in a
domain $D\subseteq{\Bbb C}$. Then
\begin{equation}\label{LI} M(f\Sigma_{\varepsilon})\geqslant\inf\limits_{\varrho\in{\rm ext\,adm}\,\Sigma_{\varepsilon}}
\int\limits_{D}\frac{\varrho^2(z)}{K_{\mu}(z)}\,dm(z)\end{equation}
for all $z_0\in\overline{D}$, where
$\varepsilon\in(0,\varepsilon_0)$, $\varepsilon_0\in(0,d_0)$,
$d_0=\sup_{z\in D}|z-z_0|$ and $\Sigma_{\varepsilon}$ denotes the
family of dashed lines consisting of all intersections of the
circles $S(z_0,r)$, $r\in(\varepsilon,\varepsilon_0)$, with $D$.}
\end{theo}

\medskip

{\it Proof.} Let $B$ be a (Borel) set of all points $z$ in $D$
where $f$ has a total differential with $J_f(z)\neq0$. It is known
that $B$ is the union of a countable collection of Borel sets
$B_l$, $l=1,2,\ldots$, such that $f_l=f|_{B_l}$ is a bi-Lipschitz
homeomorphism, see, e.g., Lemma 3.2.2 in \cite{Fe}. With no loss
of generality, we may assume that the $B_l$ are mutually disjoint.
Denote also by $B_*$ the set of all points $z\in D$ where $f$ has
a total differential with $f'(z)=0$.

Note that the set $B_0=D\setminus(B\cup B_*)$ has the Lebesgue
measure zero in ${\Bbb C}$ by the well known
Gehring--Lehto--Menchoff theorem, see \cite{GL} and \cite{Me}. Hence
by Theorem 2.11 in \cite{KR$_2$}, see also Lemma 9.1 in \cite{MRSY},
${\rm length}(\gamma\cap B_0)=0$ for a.e. paths $\gamma$ in $D$. Let
us show that ${\rm length}(f(\gamma)\cap f(B_0))=0$ for a.e. circle
$\gamma$ centered at $z_0\in\overline{D}$.

The latter follows from absolute continuity of $f$ on closed subarcs
of $\gamma\cap D$ for a.e. such circle $\gamma$. Indeed, the class
$W^{1,1}_{\rm loc}$ is invariant with respect to local
quasi-isometries and the functions in $W^{1,1}_{\rm loc}$ is
absolutely continuous on lines, see, e.g., Theorems 1.1.7 and 1.1.3
in \cite{Maz}, respectively. Applying say the transformation of
coordinates $\log(z-z_0)$, we come to the absolute continuity on
a.e. such circle $\gamma$. Fix $\gamma_0$ on which $f$ is absolutely
continuous and ${\rm length}(\gamma_0\cap B_0)=0$. Then ${\rm
length}(f(\gamma)\cap f(B_0))={\rm length}f(\gamma_0\cap B_0)$ and
for every $\varepsilon>0$ there is an open set
$\omega_{\varepsilon}$ in $\gamma_0\cap D$ such that $\gamma_0\cap
B_0\subset\omega_{\varepsilon}$ with ${\rm
length}\,\omega_{\varepsilon}<\varepsilon$, see, e.g., Theorem
III(6.6) in \cite{Sa}. The open set $\omega_{\varepsilon}$ consists
of a countable collection of open arcs $\gamma_i$ of the circle
$\gamma_0$. By the construction $\sum\limits_i{\rm
length}\,\gamma_i<\varepsilon$ and by the absolute continuity of $f$
on $\gamma_0$ the sum $\delta:=\sum\limits_i{\rm
length}\,f(\gamma_i)$ is arbitrarily small for small enough
$\varepsilon>0$. Hence ${\rm length}f(\gamma_0\cap B_0)=0$.

Thus, ${\rm length}(\gamma_*\cap f(B_0))=0$ where
$\gamma_*=f(\gamma)$ for a.e. circle $\gamma$ centered at $z_0$.
Now, let $\varrho_*\in{\rm adm}\,f(\Gamma)$ where $\Gamma$ is the
collection of all dashed lines $\gamma\cap D$ for such circles
$\gamma$ and $\varrho_*\equiv0$ outside $f(D)$. Set
$\varrho\equiv0$ outside $D$ and on $B_0$ and
$$\varrho(z)\colon=\varrho_*(f(z))\left(|f_{z}|+|f_{\bar{z}}|\right)\qquad {\rm for}\ z\in B\cup B_*\,.$$

Arguing piecewise on $B_l$, we have by Theorem 3.2.5 under $m=1$
in \cite{Fe} that $$\int\limits_{\gamma}\varrho\,ds\geqslant
\int\limits_{\gamma_*}\varrho_*\,ds_*\geqslant1\qquad {\rm for\
a.e.}\ \gamma\in\Gamma$$ because ${\rm length}(f(\gamma)\cap
f(B_0))=0$ and ${\rm length}(f(\gamma)\cap f(B_*))=0$ for a.e.
$\gamma\in\Gamma$, consequently, $\varrho\in{\rm
ext\,adm}\,\Gamma$.

On the other hand, again arguing piecewise on $B_l$, we have the
inequality $$\int\limits_{D}\frac{\varrho^2(z)}{K_{\mu}(z)}\,dm(z)
\leqslant\int\limits_{f(D)}\varrho^2_*(w)\,dm(w)$$ because
$\varrho(z)=0$ on $B_*$. Thus, we obtain (\ref{LI}).

\newpage

\begin{theo}\label{th8.4.8} {\it Let $f$ be a homeomorphic $W^{1,1}_{\rm
loc}$ solution of the Beltrami equation (\ref{eqBeltrami}) in a
domain $D\subseteq{\Bbb C}$.Then
\begin{equation}\label{eq8.4.9}
M(f\Sigma_{\varepsilon})\geqslant\int\limits_{\varepsilon}^{\varepsilon_0}
\frac{dr}{||K_{\mu}||_{1}(z_0,r)}\,,\quad\forall\
z_0\in\overline{D}\,,\ \varepsilon\in(0,\varepsilon_0)\,,\
\varepsilon_0\in(0,d_0),\end{equation} where $d_0=\sup_{z\in
D}|z-z_0|$, $\Sigma_{\varepsilon}$ denotes the family of dashed
lines consisting  of all the intersections of the circles
$S(z_0,r)$, $r\in(\varepsilon,\varepsilon_0)$, with $D$ and
\begin{equation}\label{eq8.4.11}
||K_{\mu}||_{1}(z_0,r):=\int\limits_{D(z_0,r)}K_{\mu}(z)\,
|dz|\end{equation} is the norm in $L_1$ of $K_{\mu}$ over
$D(z_0,r)=\{z\in D:|z-z_0|=r\}=D\cap S(z_0,r)$.}
\end{theo}

\medskip

{\it Proof.} Indeed, for every $\varrho\in{\rm
ext\,adm}\,\Sigma_{\varepsilon}$,
$$A_{\varrho}(r)=\int\limits_{D(z_0,r)}\varrho(z)\,|dz|\neq0\quad{\rm a.e.\ in}\ \
r\in(\varepsilon,\varepsilon_0)$$ is a measurable function in the
parameter $r$, say by the Fubini theorem. Thus, we may request the
equality $A_{\varrho}(r)\equiv1$ a.e. in
$r\in(\varepsilon,\varepsilon_0)$ and
$$\inf\limits_{\varrho\in{\rm ext\,adm}\,\Sigma_{\varepsilon}}\int\limits_{D\cap
R_{\varepsilon}}\frac{\varrho^2(z)}{K_{\mu}(z)}\,dm(z)=
\int\limits_{\varepsilon}^{\varepsilon_0}\left(\inf\limits_{\alpha\in
I(r)}\int\limits_{D(z_0,r)}\frac{\alpha^2(z)}{K_{\mu}(z)}\,|dz|\right)dr$$
where $R_{\varepsilon}=R(z_0,\varepsilon,\varepsilon_0)$ and $I(r)$
denotes the set of all measurable functions $\alpha$ on the dashed
line $D(z_0,r)=S(z_0,r)\cap D$ such that
$$\int\limits_{D(z_0,r)}\alpha(z)\,|dz|=1\,.$$
Hence Theorem \ref{th8.4.8} follows by Lemma \ref{lem8.4.1} with
$X=D(z_0,r)$, the length  as a measure $\mu$ on $D(x_0,r)$,
$\varphi=\frac{1}{K_{\mu}}|_{D(z_0,r)}$ and $p=2$.

\medskip

Later on, the following lemma will be useful, too. Here we use the
standard conventions $a/\infty=0$ for $a\neq \infty$ and
$a/0=\infty$ if $a>0$ and $a\cdot\infty=0$, see, e.g., \cite{Sa}, p.
6.

\medskip

\begin{lemma}\label{lem4cr} {\it Let $Q:{\Bbb C}\to(0,\infty)$ be a
locally integrable function. Set \begin{equation}\label{eq8.4.111}
||Q||_{1}(z_0,r)=\int\limits_{S(z_0,r)}Q(z)\, |dz|,\ \ \ z_0\in{\Bbb
C},\ r\in(0,\infty)\ ,
\end{equation}
\begin{equation}\label{eta0} I=\int\limits_{r_1}^{r_2}\frac{dr}{||Q||_1(z_0,r)}\
,\
\eta_0(r)=\frac{1}{I\cdot\|Q\|_1(z_0,r)}\,, \quad r\in(r_1,r_2),\
0<r_1<r_2<\infty\,.\end{equation} Then
\begin{equation}\label{ring}I^{-1}=\int\limits_{r_1<|z-z_0|<r_2}
Q(z)\cdot\eta_0^2\left(|z-z_0|\right)\,dm(z)\leqslant
\int\limits_{r_1<|z-z_0|<r_2}Q(z)\cdot\eta^2\left(|z-z_0|\right)\,dm(z)\end{equation}
for every measurable function $\eta:(r_1,r_2)\to[0,\infty]$ such
that
\begin{equation}\label{admr}\int\limits_{r_1}^{r_2}\eta(r)\,dr=1\, .\end{equation}

} \end{lemma}

\medskip

{\it Proof.} If $I=\infty$, then the inequality (\ref{ring}) is
obvious. Note also that $I\ne 0$ because in the contrary case
$||Q||_1(z_0,r)=\infty$ for a.e. $r\in(r_1,r_2)$ contradicting the
condition $Q\in L^1_{{\rm loc}}$. Hence we may assume further that
$0<I<\infty$. By (\ref{admr}), $\eta(r)\ne\infty$ a.e. in
$(r_1,r_2)$ and we have that $\eta(r)=\alpha(r)w(r)$ a.e. in
$(r_1,r_2)$ where
$$\alpha(r)=||Q||_1(z_0,r)\,\eta(r)\,,\quad
w(r)=\frac{1}{||Q||_1(z_0,r)}\,.$$ By the Fubini theorem in the
polar coordinates $$C\colon=\int\limits_{r_1<|z-z_0|<r_2}
Q(z)\cdot\eta^2(|z-z_0|)\,dm(z)=\int\limits_{r_1}^{r_2}\alpha^2(r)\cdot
w(r)\,dr\,.$$

By Jensen's inequality with weights $w(r)$, see, e.g., Theorem 2.6.2
in \cite{Ran} applied to the convex function $\varphi(t)=t^2$ in the
interval $\Omega=(r_1,r_2)$ and to the probability measure
$$\nu(E)=\frac{1}{I}\int\limits_E w(r)\,dr\,,\quad E\subset\Omega\,,$$
we obtain that $$\left(\dashint\alpha^2(r)w(r)\,dr\right)^{1/2}
\geqslant\dashint\alpha(r)w(r)\,dr=\frac{1}{I}$$ where we have
also used the fact that $\eta(r)=\alpha(r)\,w(r)$ satisfies
(\ref{admr}). Thus,
$$C\geqslant\left(\int\limits_{r_1}^{r_2}\frac{dr}{||Q||_1(z_0,r)}\right)^{-1}$$
and the proof is complete.

\section{On the boundary behavior of homeomorphic solutions}

\begin{theo}\label{thKPR9.1} {\it Let $D$ and $D'$ be  domains in
${\Bbb C}$ and let $f:D\to D'$ be a homeomorphic $W^{1,1}_{\rm loc}$
solution of the Beltrami equation (\ref{eqBeltrami}). Suppose that
$D$ is bounded and locally connected on $\partial D$ and $\partial
D'$ is strongly accessible. If for all $z_0\in\partial D$
\begin{equation}\label{eq8.11.2a}\int\limits_{0}^{\delta(z_0)}
\frac{dr}{||K_{\mu}||_{1}(z_0,r)}=\infty\qquad \end{equation} for
some $\delta(z_0)\in(0,d(z_0))$ where $d(z_0)=\sup_{z\in D}|z-z_0|$
and \begin{equation}\label{eq8.11.4}
||K_{\mu}||_1(z_0,r)=\int\limits_{D\cap S(z_0,r)}K_{\mu}(z)\,
|dz|\,.\end{equation} Then $f$ can be extended to $\overline{D}$ by
continuity in $\overline{\Bbb C}$.} \end{theo}

\medskip

The proof of Theorem \ref{thKPR9.1} is reduced to the following
lemma.

\medskip

\begin{lemma}\label{lem4} {\it Let $D$ and $D'$ be  domains in
${\Bbb C}$ and let $f:D\to D'$ be a homeomorphic $W^{1,1}_{\rm
loc}$ solution of the Beltrami equation (\ref{eqBeltrami}).
Suppose that the domain $D$ is bounded and locally connected at
$z_0\in\partial D$ and $\partial D'$ is strongly accessible at
least at one point of the cluster set
\begin{equation}\label{eq8.7.2} L:=C(z_0,f)=\{y\in\overline{\Bbb
C}:w=\lim\limits_{k\to\infty}f(z_k), z_k\to z_0\}\,.\end{equation}
If the condition (\ref{eq8.11.2a}) holds at $z_0$, then $f$ extends
to $z_0$ by continuity in $\overline{\Bbb C}$.}
\end{lemma}

\medskip

{\it Proof.} Note that $L\neq\varnothing$ in view of compactness of
the extended space $\overline{\Bbb C}$. By the condition $\partial
D'$ is strongly accessible at a point $\zeta_0\in L$. Let us assume
that there is one more point $z_0\in L$ and set $U=B(z_0,r_0)$ where
$0<r_0<|\zeta_0-z_0|$.

In view of local connectedness of $D$ at $z_0$, there is a sequence
of neighborhoods $V_k$ of $z_0$ with domains  $D_k=D\cap V_k$ and
$\delta(V_k)\to 0$ as $k\to\infty$. Choose in the domains
$D'_k=fD_k$ points $\zeta_k$ and $z_k$ with $|\zeta_0-\zeta_k|<r_0$
and $|\zeta_0-z_k|>r_0$, $\zeta_k\to \zeta_0$ and $z_k\to z_0$ as
$k\to\infty$. Let $C_k$ be paths connecting $\zeta_k$ and $z_k$ in
$D_k'$. Note that by the construction $\partial U\cap
C_k\neq\varnothing$. By the condition of the strong accessibility of
the point $\zeta_0$ from $D'$, there is a compactum $E\subseteq D'$
and a number $\delta>0$ such that $$M(\Delta(E,C_k;D'))\geq\delta$$
for large $k$. Without loss of generality we may assume that the
last condition holds for all $k=1,2,\ldots$. Note that $C=f^{-1}E$
is a compactum in $D'$ and hence $\varepsilon_0={\rm
dist}(z_0,C)>0$.

Let $\Gamma_{\varepsilon}$ be a family of all paths connecting the
circles $C_{\varepsilon}=\{z\in{\Bbb C}:|z-z_0|=\varepsilon\}$ and
$C_0=\{z\in\Bbb{C}:|z-z_0|=\varepsilon_0\}$ in $D$. Note that
$C_k\subset fB_{\varepsilon}$ for every fixed
$\varepsilon\in(0,\varepsilon_0)$ for large $k$ where
$B_{\varepsilon}=B(z_0,\varepsilon)$. Thus,
$M(f\Gamma_{\varepsilon})\geqslant\delta$ for all
$\varepsilon\in(0,\varepsilon_0)$. However, by \cite{He} and
\cite{Zi}, $$M(f\Gamma_{\varepsilon})\leqslant
\frac{1}{M(f\Sigma_{\varepsilon})}$$ where $\Sigma_{\varepsilon}$
is the family of all dashed lines $D(r):=\{z\in D:|z-z_0|=r\}$,
$r\in(\varepsilon,\varepsilon_0)$. Thus,
$M(f\Gamma_{\varepsilon})\to0$ as $\varepsilon\to0$ by Theorem
\ref{th8.4.8} in view of (\ref{eq8.11.2a}). The contradiction
disproves the above assumption.

\medskip

Combining Lemmas \ref{lem4cr} and \ref{lem4}, we obtain also the
following.

\medskip

\begin{lemma}\label{lem333.5.333} {\it Let $D$ and $D'$ be domains in
${\Bbb C}$, $D$ be locally connected on $\partial D$, $\partial D'$
be strongly accessible and let $Q:{\Bbb C}\to(0,\infty)$ be a
locally integrable function. Suppose that $f:D\to D'$ is a
homeomorphic $W^{1,1}_{\rm loc}$ solution of the Beltrami equation
(\ref{eqBeltrami}) such that $K_{\mu}(z)\leq Q(z)$ a.e. in $D$ and
\begin{equation}\label{omal}
\int\limits_{\varepsilon<|z-z_0|<\varepsilon_0}
Q(z)\cdot\psi^2_{z_0,\varepsilon}(|z-z_0|)\,dm(z)=o(I_{z_0}^{2}(\varepsilon))\quad{\rm
as}\quad\varepsilon\to0\ \ \forall\ z_0\in\partial D\end{equation}
for some $\varepsilon_0\in(0,\delta_0)$ where
$\delta_0=\delta(z_0)=\sup_{z\in D}|z-z_0|$ and
$\psi_{z_0,\varepsilon}(t)$ is a family of non-negative measurable
(by Lebesgue) functions on $(0,\infty)$ such that
\begin{equation}\label{eq5.3}
I_{z_0}(\varepsilon)\colon
=\int\limits_{\varepsilon}^{\varepsilon_0}
\psi_{z_0,\varepsilon}(t)\,dt<\infty\qquad\forall\
\varepsilon\in(0,\varepsilon_0)\,.\end{equation} Then $f$ can be
extended to $\overline{D}$ by continuity in $\overline{\Bbb C}$. }
\end{lemma}

\medskip

\begin{theo}\label{thBCfmo} {\it Let $D$ and $D'$ be domains in
${\Bbb C}$, $D$ be locally connected on $\partial D$, $\partial D'$
be strongly accessible and let $\mu:D\to{\Bbb C}$ be a measurable
function with $|\mu(z)|<1$ a.e. such that $K_{\mu}(z)\leqslant Q(z)$
a.e. in $D$ for a  function $Q:{\Bbb C}\to[0,\infty]$ in the class
${\rm FMO}(z_0)$ for all $z_0\in\overline{D}$. Then every
homeomorphic $W^{1,1}_{\rm loc}$ solution $f:D\to D'$ of the
Beltrami equation (\ref{eqBeltrami}) can be extended to
$\overline{D}$ by continuity in $\overline{\Bbb C}$. }
\end{theo}

\medskip

We assume that $K_{\mu}$ is extended by zero outside of $D$ in the
following consequence of Theorem \ref{thBCfmo}.

\medskip

\begin{corol}\label{corBC1} {\it Let $D$ and $D'$ be domains in
${\Bbb C}$, $D$ be locally connected on $\partial D$ and $\partial
D'$ be strongly accessible and let $\mu:D\to{\Bbb C}$ be a
measurable function with $|\mu(z)|<1$ a.e. such that
\begin{equation}\label{eqBC6*}\overline{\lim\limits_{\varepsilon\to0}}\quad
\dashint_{B(z_0,\varepsilon)}K_{\mu}(z)\,dm(z)<\infty\qquad\forall\
z_0\in\partial{D}\,.\end{equation} Then every homeomorphic
$W^{1,1}_{\rm loc}$ solution $f:D\to D'$ of the Beltrami equation
(\ref{eqBeltrami}) can be extended to $\overline{D}$ by continuity
in $\overline{\Bbb C}$.} \end{corol}

\medskip

\medskip

\begin{lemma}\label{lemKPR8.1} {\it Let $D$ and $D'$ be domains in
${\Bbb C}$, $z_1$ and $z_2$ be distinct points in $\partial D$,
$z_1\neq\infty$, and let $f:D\to D'$ be a homeomorphic
$W^{1,1}_{\rm loc}$ solution of the Beltrami equation
(\ref{eqBeltrami}). Suppose that the function $K_{\mu}$ is
integrable on the dashed lines
\begin{equation}\label{eqKPR6.1}D(r)=\{z\in D:|z-z_1|=r\}=D\cap S(z_1,r)\end{equation}
for some set $E$ of numbers $r<|z_1-z_2|$ of a positive linear measure. If $D$ is locally
connected at $z_1$ and $z_2$ and $\partial D'$ is weakly flat, then}\end{lemma}
\begin{equation}\label{eq8.10.2} C(z_1,f)\cap C(z_2,f)=\varnothing.\end{equation}

\medskip

{\it Proof.} Without loss of generality, we may assume that the
domain $D$ is bounded. Let $d=|z_1-z_2|$. Choose
$\varepsilon_0\in(0,d)$ and $\varepsilon\in(0,\varepsilon_0)$ such
that $$E_0:=\{r\in E:r\in(\varepsilon,\varepsilon_0)\}$$ has a
positive measure. The choice is possible because of a countable
subadditivity of the linear measure and because of the exhaustion
of $E$ by the sets  $$E_m:=\{r\in E:r\in(1/m,d-1/m)\}\,.$$ Note
that each of the circles $S(z_1,r)$, $r\in E_0$, separates the
points $z_1$ and $z_2$ in ${\Bbb C}$ and $D(r)$, $r\in E_0$, in
$D$. Thus, by Theorem \ref{th8.4.8} we have that
\begin{equation}\label{eq8.10.3}M(f\Sigma_{\varepsilon})>0\end{equation} where
$\Sigma_{\varepsilon}$ denotes the family of all intersections of
$D$ with the circles $$S(r)=S(z_1,r)=\{z\in{\Bbb
C}:|z-z_1|=r\}\,,\quad r\in(\varepsilon,\varepsilon_0)\,.$$

For $i=1,2$, let $C_i$ be the cluster set $C(z_i,f)$ and suppose
that $C_1\cap C_2\neq\varnothing$. Since $D$ is locally connected
at $z_1$ and $z_2$, there exist neighborhoods $U_i$ of $z_i$ such
that $W_i=D\cap U_i$, $i=1,2$ are connected and $U_1\subset
B(z_1,\varepsilon)$ and $U_2\subset{\Bbb C}\setminus
B(z_1,\varepsilon_0)$. Set
$\Gamma=\Delta(\overline{W_1},\overline{W_2};D)$. By \cite{He} and
\cite{Zi} and (\ref{eq8.10.3}) \begin{equation}\label{eq8.10.4}
M(f\Gamma)\leqslant\frac{1}{M(f\Sigma_{\varepsilon})}<\infty\,.\end{equation}
Let $\zeta_0\in C_1\cap C_2$. Without loss of generality, we may
assume that $\zeta_0\neq\infty$ because in the contrary case one
can use an additional M\"{o}bius transformation. Choose $r_0>0$
such that $S(\zeta_0,r_0)\cap fW_1\neq\varnothing$ and
$S(\zeta_0,r_0)\cap fW_2\neq\varnothing$.

By the condition $\partial D'$ is weakly flat and hence, given a
finite number $M_0>M(f\Gamma)$, there is $r_*\in(0,r_0)$ such that
$$M(\Delta(E,F;D'))\geqslant M_0$$ for all continua $E$ and $F$ in
$D'$ intersecting the circles $S(\zeta_0,r_0)$ and $S(\zeta_0,r_*)$.
However, these circles can be connected by paths $P_1$ and $P_2$ in
the domains $fW_1$ and $fW_2$, respectively, and for those paths
$$M_0\leqslant M(\Delta(P_1,P_2;D'))\leqslant M(f\Gamma)\,.$$

The contradiction disproves the above assumption that $C_1\cap
C_2\neq\varnothing$. The proof is complete.

\medskip

As an immediate consequence of Lemma \ref{lemKPR8.1}, we have the
following statement.

\medskip

\begin{theo}\label{thKPR8.2} {\it Let $D$ and $D'$ be domains in
${\Bbb C}$, $D$ locally connected on $\partial D$ and $\partial
D'$ weakly flat. If $f:D\to D'$ is a homeomorphic $W^{1,1}_{\rm
loc}$ solution of the Beltrami equation (\ref{eqBeltrami}) with
$K_{\mu}\in L^1(D)$, then $f^{-1}$ has an extension to
$\overline{D'}$ by continuity in $\overline{\Bbb C}$.} \end{theo}

\medskip

{\it Proof.} By the Fubini theorem with notations from Lemma
\ref{lemKPR8.1}, the set \begin{equation}\label{eqKPR6.2a}
E=\{r\in(0,d):K_{\mu}|_{D(r)}\in L^{1}(D(r))\}\end{equation} has a
positive linear measure because $K_{\mu}\in L^1(D)$.

\medskip

\begin{rem}\label{rmkKPR6.1a} It is clear that it is even
sufficient to assume in Theorem \ref{thKPR8.2} that $K_{\mu}$ is
integrable only in a neighborhood of $\partial D$. \end{rem}

\medskip

Moreover, by Lemma \ref{lemKPR8.1} we obtain also the following
conclusion.

\medskip

\begin{theo}\label{thKPR8.3} {\it Let $D$ and $D'$ be domains in
${\Bbb C}$, $D$ bounded and locally connected on $\partial D$ and
$\partial D'$ weakly flat, and let $f:D\to D'$ be a homeomorphic
$W^{1,1}_{\rm loc}$ solution of the Beltrami equation
(\ref{eqBeltrami}) with the coefficient $\mu$ such that the
condition (\ref{eq8.11.2a}) holds for all $z_0\in\partial D$. Then
there is an extension of $f^{-1}$ to $\overline{D'}$ by continuity
in $\overline{\Bbb C}$.} \end{theo}

\medskip

Combining Theorems \ref{thKPR9.1} and \ref{thKPR8.3}, we obtain the following
statement.

\medskip

\begin{theo}\label{thKPR8.3p} {\it Let $D$ and $D'$ be domains in
${\Bbb C}$ and let $f:D\to D'$ be a homeomorphic $W^{1,1}_{\rm
loc}$ solution of the Beltrami equation (\ref{eqBeltrami}).
Suppose that $D$ is bounded and locally connected on $\partial D$
and $\partial D'$ is weakly flat. If the condition
(\ref{eq8.11.2a}) holds for all $z_0\in\partial D$, then $f$ has a
homeomorphic extension $\overline{f}:\overline{D}\to\overline{D'}$
by continuity in $\overline{\Bbb{C}}$.} \end{theo}

\medskip

In particular, as a consequence of Theorem \ref{thKPR8.3p} we obtain
in the plane the following generalization of the well-known
Gehring--Martio theorem on a homeomorphic extension to the boundary
of quasiconformal mappings between QED domains in ${\Bbb R}^n$,
$n\geq 2$, cf. \cite{GM}, see also \cite{MV}.

\medskip

\begin{corol}\label{thKPR9.2} {\it Let $D$ and $D'$ be bounded domains
with weakly flat boundaries in ${\Bbb C}$ and let $f:D\to D'$ be a
homeomorphic $W^{1,1}_{\rm loc}$ solution of the Beltrami equation
(\ref{eqBeltrami}). If the condition (\ref{eq8.11.2a}) holds at
every point $z_0\in\partial D$, then $f$ has a homeomorphic
extension to $\overline{D}$ by continuity in $\overline{\Bbb C}$.}
\end{corol}

\medskip

By Theorem \ref{thKPR8.3p} and Lemma \ref{lem4cr} we have also the
following.

\medskip

\begin{lemma}\label{lem13.5.333} {\it Let $D$ and $D'$ be domains in
${\Bbb C}$, $D$ be locally connected on $\partial D$, $\partial D'$
be weakly flat and let $Q:{\Bbb C}\to(0,\infty)$ be a locally
integrable function. Suppose that $f:D\to D'$ is a homeomorphic
$W^{1,1}_{\rm loc}$ solution of the Beltrami equation
(\ref{eqBeltrami}) such that $K_{\mu}(z)\leq Q(z)$ a.e. in $D$ and
\begin{equation}\label{omal}
\int\limits_{\varepsilon<|z-z_0|<\varepsilon_0}
Q(z)\cdot\psi^2_{z_0,\varepsilon}(|z-z_0|)\,dm(z)=o(I_{z_0}^{2}(\varepsilon))\quad{\rm
as}\quad\varepsilon\to0\ \ \forall\ z_0\in\partial D\end{equation}
for some $\varepsilon_0\in(0,\delta_0)$ where
$\delta_0=\delta(z_0)=\sup_{z\in D}|z-z_0|$ and
$\psi_{z_0,\varepsilon}(t)$ is a family of non-negative measurable
(by Lebesgue) functions on $(0,\infty)$ such that
\begin{equation}\label{eq5.3}
I_{z_0}(\varepsilon)\colon
=\int\limits_{\varepsilon}^{\varepsilon_0}
\psi_{z_0,\varepsilon}(t)\,dt<\infty\qquad\forall\
\varepsilon\in(0,\varepsilon_0)\,.\end{equation} Then $f$ has a
homeomorphic extension $\overline{f}:\overline{D}\to\overline{D'}$
by continuity in $\overline{\Bbb C}$.} \end{lemma}

\medskip

\begin{theo}\label{thBHfmo} {\it Let $D$ and $D'$ be domains in
${\Bbb C}$, $D$ be locally connected on $\partial D$, $\partial D'$
be weakly flat and let $\mu:D\to{\Bbb C}$ be a measurable function
with $|\mu(z)|<1$ a.e. such that $K_{\mu}(z)\leqslant Q(z)$ a.e. in
$D$ for a  function $Q:{\Bbb C}\to[0,\infty]$ in the class ${\rm
FMO}(z_0)$ for all $z_0\in\overline{D}$. Then every homeomorphic
$W^{1,1}_{\rm loc}$ solution $f:D\to D'$ of the Beltrami equation
(\ref{eqBeltrami}) has a homeomorphic extension
$\overline{f}:\overline{D}\to\overline{D'}$ by continuity in
$\overline{\Bbb C}$.}
\end{theo}

\medskip

We assume that $K_{\mu}$ is extended by zero outside of $D$ in the
following consequence of Theorem \ref{thBHfmo}.

\medskip

\begin{corol}\label{corBH1} {\it Let $D$ and $D'$ be domains in
${\Bbb C}$, $D$ be locally connected on $\partial D$ and $\partial
D'$ be weakly flat and let $\mu:D\to{\Bbb C}$ be a measurable
function with $|\mu(z)|<1$ a.e. such that
\begin{equation}\label{eqBH6*}\overline{\lim\limits_{\varepsilon\to0}}\quad
\dashint_{B(z_0,\varepsilon)}K_{\mu}(z)\,dm(z)<\infty\qquad\forall\
z_0\in\partial{D}\,.\end{equation} Then every homeomorphic
$W^{1,1}_{\rm loc}$ solution $f:D\to D'$ of the Beltrami equation
(\ref{eqBeltrami}) has a homeomorphic extension
$\overline{f}:\overline{D}\to\overline{D'}$ by continuity in
$\overline{\Bbb C}$.} \end{corol}

\medskip

\section{On regular solutions for the Dirichlet problem in the Jordan domains}

If $\varphi(\zeta)\not\equiv{\rm const}$, then the {\bf regular
solution} of such a problem is a continuous, discrete and open
mapping $f:D\to{\Bbb C}$ of the Sobolev class $W_{\rm loc}^{1,1}$
with its Jacobian $J_f(z)=|f_z|^2-|f_{\bar z}|^2\neq0$ a.e.
satisfying (\ref{eqBeltrami}) a.e. and the condition
(\ref{eqGrUsl}). The regular solution of the Dirichlet problem
(\ref{eqGrUsl}) with $\varphi(\zeta)\equiv c$, $\zeta\in\partial
D$, for the Beltrami equation (\ref{eqBeltrami}) is the function
$f(z)\equiv c$, $z\in D$.

Recall that a mapping $f:D\to{\Bbb C}$ is called {\bf discrete} if
the preimage  $f^{-1}(y)$ consists of isolated points for every
$y\in{\Bbb C}$, and {\bf open} if $f$ maps every open set
$U\subseteq D$ onto an open set in ${\Bbb C}$.

In this section, we prove that a regular solution of the Dirichlet
problem (\ref{eqGrUsl}) exists for every continuous function
$\varphi:\partial D\to{\Bbb R}$ for wide classes of the degenerate
Beltrami equations (\ref{eqBeltrami}) in an arbitrary  Jordan
domain $D$.

\medskip

\begin{lemma}\label{lemDIR9} {\it Let $D$ be a Jordan domain in
${\Bbb C}$ and $\mu:D\to{\Bbb C}$ be a measurable function with
$|\mu(z)|<1$ a.e. If $K_{\mu}(z)\leq Q(z)$ a.e. in $D$ where
$Q:{\Bbb C}\to(0,\infty)$ is locally integrable and satisfies the
condition (\ref{omal}) for all $z_0\in\overline{D}$, then the
Beltrami equation (\ref{eqBeltrami}) has a regular solution $f$ of
the Dirichlet problem (\ref{eqGrUsl}) for each continuous function
$\varphi:\partial D\to{\Bbb R}$.}\end{lemma}

\medskip

{\it Proof}. Let $F$ be a regular homeomorphic solution of the
Beltrami equation (\ref{eqBeltrami}) of the class $W^{1,1}_{\rm
loc}$ that exists by Lemma 4.1 in \cite{RSY$_6$}. Note that
$\overline{\Bbb C}\setminus D^*$, where $D^*=F(D)$, cannot consist
of the single point $\infty$ because in the contrary case $\partial
D^*$ would be weakly flat. But then by Lemma \ref{lem13.5.333} $F$
should have a homeomorphic extension to $\overline{D}$ that is
impossible because $\partial D$ is not a singleton. Moreover, the
domain $D^*$ is simply connected, see, e.g., either Lemma 5.3 in
\cite{IR} or Lemma 6.5 in \cite{MRSY}. Thus, by the Riemann theorem,
see, e.g., Theorem II.2.1 in \cite{Gol}, $D^*$ can be mapped by a
conformal mapping $R$ onto the unit disk ${\Bbb D}$. The mapping
$g=R\circ F$ is also a regular homeomorphic solution of the Beltrami
equation of the class $W_{\rm loc}^{1,1}$ that maps $D$ onto ${\Bbb
D}$. Furthermore, by Lemma \ref{lem13.5.333} $g$ admits a
homeomorphic extension $g_*:D\to\overline{\Bbb D}$ because ${\Bbb
D}$ has a weakly flat boundary and the Jordan domain $D$ is locally
connected on its boundary.

Let us find a solution of the Dirichlet problem (\ref{eqGrUsl}) in
the form $f=h\circ g$ where $h$ is an analytic function in ${\Bbb D}$
with the boundary condition $$\lim\limits_{z\to\zeta}{\rm Re}\,h(z)=\varphi(g_{*}^{-1}(\zeta))
\qquad\forall\ \zeta\in\partial{\Bbb D}\,.$$

By the Schwarz formula (see, e.g., Section 8, Chapter III, Part 3
in \cite{HuCo}), the analytic function $h$ with ${\rm Im}\,h(0)=0$
can be calculated in ${\Bbb D}$ through its real part on the
boundary: \begin{equation}\label{eqDIR4*} h(z)=\frac{1}{2\pi
i}\int\limits_{|\zeta|=1}{\rm Re}\,\varphi\circ
g_{*}^{-1}(\zeta)\cdot\frac{\zeta+z}{\zeta-z}\cdot\frac{d\zeta}{\zeta}\,.\end{equation}

We see that the function $f=h\circ g$ is the desired regular
solution of the Dirichlet problem (\ref{eqGrUsl}) for the Beltrami
equation (\ref{eqBeltrami}).

\medskip

Choosing in Lemma \ref{lemDIR9} $\psi(t)=1/\left(t\,
\log\left(1/t\right)\right)$, we obtain by Lemma \ref{lem13.4.2} the
following result.

\medskip

\begin{theo}\label{thDIR2fmo} {\it Let $D$ be a Jordan domain and $\mu:D\to{\Bbb C}$
be a measurable function with $|\mu(z)|<1$ a.e. such that
$K_{\mu}(z)\leqslant Q(z)$ a.e. in $D$ for a  function $Q:{\Bbb
C}\to[0,\infty]$ in ${\rm FMO}(\overline{D})$. Then the Beltrami
equation (\ref{eqBeltrami}) has a regular solution of the Dirichlet
problem (\ref{eqGrUsl}) for each continuous function
$\varphi:\partial D\to{\Bbb R}$.} \end{theo}

\medskip

\begin{corol}\label{corDIR1000re} {\it In particular, the conclusion of
Theorem \ref{thDIR2fmo} holds if every point $z_0\in\overline{D}$ is
the Lebesgue point of a locally integrable function $Q:{\Bbb
C}\to[0,\infty]$ such that $K_{\mu}(z)\leqslant Q(z)$ a.e. in $D$.}
\end{corol}

\medskip

Further we assume that $K_{\mu}$ is extended by zero outside of $D$.

\medskip

\begin{corol}\label{corDIR1} {\it Let $D$ be a Jordan domain
and $\mu:D\to{\Bbb C}$ be a measurable function with $|\mu(z)|<1$
a.e. such that
\begin{equation}\label{eqDIR6*}\overline{\lim\limits_{\varepsilon\to0}}\quad
\dashint_{B(z_0,\varepsilon)}K_{\mu}(z)\,dm(z)<\infty\qquad\forall\
z_0\in\overline{D}\,.\end{equation} Then the Beltrami equation
(\ref{eqBeltrami}) has a regular solution of the Dirichlet problem
(\ref{eqGrUsl}) for each continuous function $\varphi:\partial
D\to{\Bbb R}$.} \end{corol}

\medskip

The following statement is proved similarly to Lemma \ref{lemDIR9}
on the basis of Theorem \ref{thKPR8.3p} instead of Lemma
\ref{lem13.5.333}.

\medskip

\begin{theo}\label{thDIR2io} {\it Let $D$ be a Jordan domain in $\Bbb{C}$ and $\mu:D\to{\Bbb C}$
be a measurable function with $|\mu(z)|<1$ a.e. If $K_{\mu}\in
L^1_{\rm loc}(D)$ and satisfies the condition
\begin{equation}\label{eq8.11.2}\int\limits_{0}^{\delta(z_0)}
\frac{dr}{||K_{\mu}||_{1}(z_0,r)}=\infty\qquad\forall\
z_0\in\overline{D}\end{equation} for some $\delta(z_0)\in(0,d(z_0))$
where $d(z_0)=\sup\limits_{z\in D}|z-z_0|$ and
\begin{equation}\label{eq8.11.4}
||K_{\mu}||_{1}(z_0,r)=\int\limits_{D\cap S(z_0,r)}K_{\mu}(z)\,
|dz|\,,\end{equation} at each point $z_0\in\overline{D}$, then the
Beltrami equation (\ref{eqBeltrami}) has a regular solution $f$ of
the Dirichlet problem (\ref{eqGrUsl}) for each continuous function
$\varphi:\partial D\to{\Bbb R}$.}\end{theo}

\medskip

\begin{corol}\label{corDIR2t} {\it Let $D$ be a Jordan domain
and $\mu:D\to{\Bbb C}$ be a measurable function such that
\begin{equation}\label{eqDIR61}k_{z_{0}}(\varepsilon)=O\left(\log\frac{1}{\varepsilon}\right)
\qquad\forall\ z_0\in\overline{D}\end{equation} as
$\varepsilon\to0$, where $k_{z_0}(\varepsilon)$ is the average of
the function $K_{\mu}(z)$ over $S(z_{0},\varepsilon)$. Then the
Beltrami equation (\ref{eqBeltrami}) has a regular solution of the
Dirichlet problem (\ref{eqGrUsl}) for each continuous function
$\varphi:\partial D\to{\Bbb R}$.} \end{corol}

\medskip

\begin{rem}\rm\label{rem2} {\rm In particular, the conclusion of Corollary \ref{corDIR2t} holds if
\begin{equation}\label{eqDIR6**} K_{\mu}(z)=O\left(\log\frac{1}{|z-z_0|}\right)\qquad{\rm
as}\quad z\to z_0\quad\forall\
z_0\in\overline{D}\,.\end{equation}}\end{rem}

\medskip

Finally, combining Theorems \ref{th5.555} and \ref{thDIR2io} we
obtain the following.

\medskip

\begin{theo}\label{thKR4.1} {\it Let $D$ be a Jordan domain and
$\mu:D\to{\Bbb D}$ be a measurable function such that
\begin{equation}\label{eqKR4.1sh}\int\limits_{D}\Phi\left(K_{\mu}(z)\right)\,dm(z)<\infty\end{equation}
for a convex non-decreasing function
$\Phi:[0,\infty]\to[0,\infty]$. If
\begin{equation}\label{eqKR4.2ras}
\int\limits_{\delta}^{\infty}\frac{d\tau}{\tau\Phi^{-1}(\tau)}=\infty\end{equation}
for some $\delta>\Phi(0)$. Then the Beltrami equation
(\ref{eqBeltrami}) has a regular solution of the Dirichlet problem
(\ref{eqGrUsl}) for each continuous function $\varphi:\partial
D\to{\Bbb R}$.} \end{theo}

\medskip

\begin{rem}\label{remeq333F} By the Stoilow theorem, see, e.g., \cite{Sto}, a regular solution $f$
of the Dirichlet problem (\ref{eqGrUsl}) for the Beltrami equation
(\ref{eqBeltrami}) with $K_{\mu}\in L^1_{\rm loc}(D)$ can be
represented in the form $f=h\circ F$ where $h$ is an analytic
function and $F$ is a homeomorphic regular solution of
(\ref{eqBeltrami}) in the class $W_{\rm loc}^{1,1}$. Thus, by
Theorem 5.1 in \cite{RSY13} the condition (\ref{eqKR4.2ras}) is
not only sufficient but also necessary to have a regular solution
of the Dirichlet problem (\ref{eqGrUsl}) for an arbitrary Beltrami
equation (\ref{eqBeltrami}) with the integral constraints
(\ref{eqKR4.1sh}) for any non-constant continuous  function
$\varphi:\partial D\to\Bbb{R}$.

Setting $H(t)=\log\Phi(t)$, note that by Theorem 2.1 in \cite{RSY}
the condition \ref{eqKR4.2ras} is equivalent to each of the
conditions
\begin{equation}\label{eq333Frer}\int\limits_{\Delta}^{\infty}H'(t)\,\frac{dt}{t}=\infty,
\end{equation} \begin{equation}\label{eq333F}\int\limits_{\Delta}^{\infty}
\frac{dH(t)}{t}=\infty\,,\end{equation} and (\ref{eq333F})
\begin{equation}\label{eq333B}
\int\limits_{\Delta}^{\infty}H(t)\,\frac{dt}{t^2}=\infty\,\end{equation}
for some $\Delta>0$, and \begin{equation}\label{eq333C}
\int\limits_{0}^{\delta}H\left(\frac{1}{t}\right)\,{dt}=\infty\end{equation}
for some $\delta>0$, \begin{equation}\label{eq333D}
\int\limits_{\Delta_*}^{\infty}\frac{d\eta}{H^{-1}(\eta)}=\infty\end{equation}
for some $\Delta_*>H(+0)$. Here, the integral in (\ref{eq333F}) is
understood as the Lebesgue--Stieltjes integral and the integrals
in (\ref{eqKR4.2ras}) and (\ref{eq333B})--(\ref{eq333D}) as the
ordinary Lebesgue integrals. \end{rem}

\medskip

\begin{corol}\label{corDIR1000} {\it In particular, the conclusion of
Theorem \ref{thKR4.1} holds if, for some $\alpha>0$,
\begin{equation}\label{eqKR4.1}\int\limits_{D}e^{\alpha K_{\mu}(z)}\,dm(z)<\infty\,.\end{equation}}
\end{corol}

\section{On pseudoregular solutions in multiply connected domains}

As it was first noted by Bojarski, see, e.g., section 6 of Chapter
4 in \cite{Vekua}, in the case of multiply connected domains the
Dirichlet problem for the Beltrami equation, generally speaking,
has no solutions in the class of continuous (simply-valued)
functions. Hence it is arose the question: whether the existence
of solutions for the Dirichlet problem can be obtained in a wider
class for the case? It is turned out to be that this is possible
in the class of functions having a certain number of poles at
prescribed points in $D$. More precisely, for
$\varphi(\zeta)\not\equiv{\rm const}$, a {\bf pseudoregular
solution} of the problem is a continuous (in $\overline{\Bbb
C}={\Bbb C}\cup\{\infty\}$) discrete open mapping
$f:D\to\overline{\Bbb C}$ in the class $W^{1,1}_{\rm loc}$
(outside of these poles) with the Jacobian
$J_{f}(z)=|f_z|^2-|f_{\bar z}|^2\neq0$ a.e. satisfying
(\ref{eqBeltrami}) a.e. and the condition (\ref{eqGrUsl}).

\medskip

\begin{lemma}\label{lem13.5.333ps} {\it Let $D$ be a bounded domain in
${\Bbb C}$ whose boundary consists of $n\geqslant2$ mutually
disjoint Jordan curves and let $\mu:D\to{\Bbb C}$ be a measurable
function such that $|\mu(z)|<1$ a.e. Suppose that $K_{\mu}(z)\leq
Q(z)$ a.e. in $D$ where $Q:{\Bbb C}\to(0,\infty)$ is locally
integrable and satisfies the condition (\ref{omal}) for all
$z_0\in\overline{D}$. Then the Beltrami equation (\ref{eqBeltrami})
has a pseudoregular solution of the Dirichlet problem
(\ref{eqGrUsl}) for each continuous function $\varphi:\partial
D\to{\Bbb R}$, $\varphi(\zeta)\not\equiv{\rm const}$, with poles at
$n$ prescribed points $z_i\in D$, $i=1,\ldots,n$.} \end{lemma}

\medskip

{\it Proof.} Let $F$ be a regular homeomorphic solution of the
Beltrami equation (\ref{eqBeltrami}) of the class $W^{1,1}_{\rm loc}$ that
exists by Lemma 4.1 in \cite{RSY$_6$}. Consider $D^*=f(D)$.
Note that $\partial D^*$ has $n$ connected components $\Gamma_i$, $i=1,\ldots,n$
that correspond in the natural way to connected components of $\partial D$,
the Jordan curves $\gamma_i$, see, e.g., either Lemma 5.3 in \cite{IR}
or Lemma 6.5 in \cite{MRSY}.

Thus, by Theorem V.6.2 in \cite{Gol} the domain $D_*$ can be mapped
by a conformal map $R$ onto a circular domain ${\Bbb D}_*$ whose
boundary consists of $n$ circles or points, i.e. ${\Bbb D}_*$ has a
weakly flat boundary. Note that the mapping $g:=R\circ F$ is a
regular homeomorphic solution of the Beltrami equation in  the class
$W_{\rm loc}^{1,1}$ admitting a homeomorphic extension
$g_*:\overline{D}\to\overline{\Bbb D_*}$ by Lemma \ref{lem13.5.333}.

Let us find a solution of the Dirichlet problem (\ref{eqGrUsl}) in
the form $f=h\circ g$ where $h$ is a meromorphic function with $n$
poles at the prescribed points $w_i=g(z_i)$, $i=1,\ldots,n$ in
${\Bbb D}_*$ with the boundary condition
$$\lim\limits_{w\to\zeta}{\rm Re}\,h(w)=\varphi(g_{*}^{-1}(\zeta))
\qquad\forall\ \zeta\in\partial{\Bbb D}_*$$ Such a function $h$
exists by theorem 4.14 in \cite{Vekua}.

We see that the function $f=h\circ g$ is the desired pseudoregular
solution of the Dirichlet problem (\ref{eqGrUsl}) for the Beltrami
equation (\ref{eqBeltrami}) with $n$ poles  just at these
prescribed points $z_i$, $i=1,\ldots,n$.

\medskip

Arguing, similarly to the last section, by the special choice of
the functional parameter $\psi$ in Lemma \ref{lem13.5.333ps}, we
obtain the following result.

\medskip

\begin{theo}\label{thDIR2} {\it Let $D$ be a bounded domain in
${\Bbb C}$ whose boundary consists of $n\geqslant2$ mutually
disjoint Jordan curves and $\mu:D\to{\Bbb C}$ be a measurable
function such that $|\mu(z)|<1$ a.e., and $K_{\mu}(z)\leqslant
Q(z)$ a.e. in $\overline{D}$ for a  function $Q:{\Bbb
C}\to[0,\infty]$ in the class ${\rm FMO}(\overline{D})$. Then the
Beltrami equation (\ref{eqBeltrami}) has a pseudoregular solution
of the Dirichlet problem (\ref{eqGrUsl}) for every continuous
function $\varphi:\partial D\to{\Bbb R}$,
$\varphi(\zeta)\not\equiv{\rm const}$, with poles at $n$
prescribed points in $D$.} \end{theo}

\medskip

\begin{corol}\label{corDIR1500re} {\it In particular, the conclusion of
Theorem \ref{thDIR2} holds if every point $z_0\in\overline{D}$ is
the Lebesgue point of a locally integrable function $Q:{\Bbb
C}\to[0,\infty]$ such that $K_{\mu}(z)\leqslant Q(z)$ a.e. in
$D$.}\end{corol}

\medskip

\begin{corol}\label{corDIR1000} {\it In particular, the conclusion of
Theorem \ref{thDIR2} holds if \begin{equation}\label{eqDIR6*}
\overline{\lim\limits_{\varepsilon\to0}}\ \ \
\dashint_{B(z_0,\varepsilon)}K_{\mu}(z)\,dm(z)<\infty\qquad\forall\
z_0\in\overline{D}\,.\end{equation}} \end{corol}

\medskip

As above, here we assume that $K_{\mu}$ is extended by zero outside
of $D$.

\medskip

\begin{theo}\label{thDIR2i} {\it Let $D$  be a bounded  domain in $\Bbb{C}$ whose
boundary consists of $n\geqslant2$ mutually disjoint Jordan curves
and $\mu:D\to{\Bbb C}$ be a measurable function with $|\mu(z)|<1$
a.e. If $K_{\mu}\in L^1_{\rm loc}(D)$ and satisfies the condition
\begin{equation}\label{eq8.11.2}\int\limits_{0}^{\delta(z_0)}
\frac{dr}{||K_{\mu}||_{1}(z_0,r)}=\infty\qquad\forall\
z_0\in\overline{D}\end{equation} for some
$\delta(z_0)\in(0,d(z_0))$ where $d(z_0)=\sup\limits_{z\in
D}|z-z_0|$ and
\begin{equation}\label{eq8.11.4}
||K_{\mu}||_{1}(z_0,r)=\int\limits_{D\cap S(z_0,r)}K_{\mu}\,
|dz|\,.\end{equation} Then the Beltrami equation (\ref{eqBeltrami})
has a pseudoregular solution of the Dirichlet problem
(\ref{eqGrUsl}) for every continuous function $\varphi:\partial
D\to{\Bbb R}$, $\varphi(\zeta)\not\equiv{\rm const}$, with poles at
$n$ prescribed points in $D$.} \end{theo}

\medskip

\begin{corol}\label{corDIR2t} {\it Let $D$ be a bounded domain in
${\Bbb C}$ whose boundary consists of $n\geqslant2$ mutually
disjoint Jordan curves and $\mu:D\to{\Bbb C}$ be a measurable
function such that \begin{equation}\label{eqDIR61fa}
k_{z_0}(\varepsilon)=O\left(\log\frac{1}{\varepsilon}\right)\quad
\text{as}\quad\varepsilon\to0\quad\forall\
z_0\in\overline{D}\,,\end{equation} where $k_{z_0}(\varepsilon)$
is the average of the function $K_{\mu}(z)$ over the circle
$\{z\in{\Bbb C}:|z-z_0|=\varepsilon\}$. Then the Beltrami equation
(\ref{eqBeltrami}) has a pseudoregular solution of the Dirichlet
problem (\ref{eqGrUsl}) for every continuous function
$\varphi:\partial D\to{\Bbb R}$, $\varphi(\zeta)\not\equiv{\rm
const}$, with poles at $n$ prescribed points in $D$.}\end{corol}

\medskip

\begin{rem}\label{rem2tyyr} In particular, the conclusion of Corollary \ref{corDIR2t} holds if
\begin{equation}\label{eqDIR6w**}K_{\mu}(z)=O\left(\log\frac{1}{|z-z_0|}\right)\qquad{\rm
as}\quad z\to z_0\quad\forall\ z_0\in\overline{D}\,.\end{equation}
\end{rem}

\medskip

\begin{theo}\label{thq13.5.333} {\it Let $D$ be a bounded domain in
${\Bbb C}$ whose boundary consists of $n\geqslant2$ mutually
disjoint Jordan curves and let $\mu:D\to{\Bbb C}$ be a measurable
function with $|\mu(z)|<1$ a.e. such that
\begin{equation}\label{eqKR4.1}\int\limits_{D}\Phi\left(K_{\mu}(z)\right)\,dm(z)<\infty\end{equation}
for a convex non-decreasing function $\Phi:[0,\infty]\to[0,\infty]$. If
\begin{equation}\label{eqKR4.2}\int\limits_{\delta}^{\infty}\frac{d\tau}{\tau\Phi^{-1}(\tau)}=\infty\end{equation}
for some $\delta>\Phi(0)$. Then the Beltrami equation
(\ref{eqBeltrami}) has a pseudoregular solution of the Dirichlet
problem (\ref{eqGrUsl}) for each continuous function
$\varphi:\partial D\to{\Bbb R}$, $\varphi(\zeta)\not\equiv{\rm
const}$, with poles at $n$ prescribed inner points in $D$.}
\end{theo}

\medskip

\begin{corol}\label{corDIR1000} {\it In particular, the conclusion of
Theorem \ref{thq13.5.333} holds if, for some $\alpha>0$,
\begin{equation}\label{eqKR4.1}\int\limits_{D}e^{\alpha K_{\mu}(z)}\,dm(z)<\infty\,.\end{equation}}\end{corol}

\section{On multi-valued solutions in multiply connected domains}

In  finitely connected domains $D$ in $\Bbb{C}$, in addition to
pseudoregular of solutions, the Dirichlet problem (\ref{eqGrUsl})
for the Beltrami equation (\ref{eqBeltrami}) admits multi-valued
solutions in the spirit of the theory of multi-valued analytic
functions. We say that a discrete open mapping
$f:B(z_0,\varepsilon_0)\to{\Bbb C}$, where
$B(z_0,\varepsilon_0)\subseteq D$, is a {\bf local regular solution
of the equation} (\ref{eqBeltrami}) if $f\in W_{\rm loc}^{1,1}$,
$J_f(z)\neq0$ and $f$ satisfies (\ref{eqBeltrami}) a.e. in
$B(z_0,\varepsilon_0)$.

\medskip

The local regular solutions $f:B(z_0,\varepsilon_0)\to{\Bbb C}$
and $f_*:B(z_*,\varepsilon_*)\to{\Bbb C}$ of the equation
(\ref{eqBeltrami}) will be called extension of each to other if
there is a finite chain of such solutions
$f_i:B(z_i,\varepsilon_i)\to\Bbb{C}$, $i=1,\ldots,m$, that
$f_1=f_0$, $f_m=f_*$ and $f_i(z)\equiv f_{i+1}(z)$ for $z\in
E_i:=B(z_i,\varepsilon_i)\cap
B(z_{i+1},\varepsilon_{i+1})\neq\emptyset$, $i=1,\ldots,m-1$. A
collection of local regular solutions
$f_j:B(z_j,\varepsilon_j)\to{\Bbb C}$, $j\in J$, will be called a
{\bf multi-valued solution} of the equation (\ref{eqBeltrami}) in
$D$ if the disks $B(z_j,\varepsilon_j)$ cover the whole domain $D$
and $f_j$ are extensions of each to other through the collection.
A multi-valued solution of the equation (\ref{eqBeltrami}) will be
called a {\bf multi-valued solution of the Dirichlet problem}
(\ref{eqGrUsl}) if $u(z)={\rm Re}\,f(z)={\rm Re}\,f_{j}(z)$, $z\in
B(z_j,\varepsilon_j)$, $j\in J$, is a simply-valued function in
$D$ satisfying the condition
$\lim\limits_{z\in\zeta}u(z)=\varphi(\zeta)$ for all $\zeta\in D$.

\medskip

\begin{lemma}\label{lem13.5.3334} {\it Let $D$ be a bounded domain in
${\Bbb C}$ whose boundary consists of $n\geqslant2$ mutually
disjoint Jordan curves and let $\mu:D\to{\Bbb C}$ be a measurable
function such that $|\mu(z)|<1$ a.e. Suppose that $K_{\mu}(z)\leq
Q(z)$ a.e. in $D$ where $Q:{\Bbb C}\to(0,\infty)$ is locally
integrable and satisfies (\ref{omal}) for all $z_0\in\overline{D}$.
Then the Beltrami equation (\ref{eqBeltrami}) has a multi-valued
solutions of the Dirichlet problem (\ref{eqGrUsl}) for each
continuous function $\varphi:\partial D\to{\Bbb R}$.} \end{lemma}

\medskip

{\it Proof}. Let $F$ be a regular homeomorphic solution of the
Beltrami equation (\ref{eqBeltrami}) in the class $W_{\rm
loc}^{1,1}$ that exists by Lemma 4.1 in \cite{RSY$_6$}. As it was
showed under the proof of Lemma \ref{lem13.5.333}, we may assume
that $D_*:=F(D)$ is a circular domain and that $F$ can be extended
to a homeomorphism $F_*:\overline{D}\to\overline{D_*}$. Let
$u:D_*\to\Bbb{R}$ be a harmonic function such that
$$\lim\limits_{w\to\zeta}u(w)=\varphi (F_{*}^{-1}(\zeta))\qquad\forall\ \zeta\in\partial D^*$$
whose existence is well-known, see, e.g., Section 3 of Chapter VI
in \cite{Gol}.

Let $z_0\in D$, $B_0:=B(z_0,\varepsilon_0)\subseteq D$ for some
$\varepsilon_0>0$. Then the domain $D_0=F(B_0)$ is simply connected
and hence there is a harmonic function $v(w)$ such that
$h(w)=u(w)+iv(w)$ is a holomorphic function which is unique up to an
additive constant, see, e.g., Theorem 1 in Section 7 of Chapter III,
Part 3 in \cite{HuCo}. Note that $f_0:=h\circ F|_{B_0}$ is a local
regular solution of the Beltrami equation (\ref{eqBeltrami}). Note
that the function $h$ can be extended to a multi-valued analytic
function $H$ in the domain $D_*$ and, thus, $H\circ F$ gives the
desired multi-valued solutions of the Dirichlet problem
(\ref{eqGrUsl}) for the Beltrami equation (\ref{eqBeltrami}).

\medskip

In particular, by Lemmas \ref{lem13.4.2} and \ref{lemDIR9} we
obtain the following.

\medskip

\begin{theo}\label{12thDIR2fmo} {\it Let $D$ be a bounded domain in
${\Bbb C}$ whose boundary consists of $n\geqslant2$ mutually
disjoint Jordan curves and $\mu:D\to{\Bbb C}$ be a measurable
function with $|\mu(z)|<1$ a.e. and such that $K_{\mu}(z)\leqslant
Q(z)$ a.e. in $D$ for a function $Q:\Bbb{C}\to[0,\infty]$ in ${\rm
FMO}(\overline{D})$. Then the Beltrami equation (\ref{eqBeltrami})
has a multi-valued solutions of the Dirichlet problem
(\ref{eqGrUsl}) for each continuous function $\varphi:\partial
D\to{\Bbb R}$.}\end{theo}

\medskip

\begin{corol}\label{12corDIR1500re} {\it In particular, the conclusion of
Theorem \ref{12thDIR2fmo} holds if every point
$z_0\in\overline{D}$ is Lebesgue point of a locally integrable
function $Q:{\Bbb C}\to[0,\infty]$ such that $K_{\mu}(z)\leqslant
Q(z)$ a.e. in $D$.}
\end{corol}

\medskip

We assume that $K_{\mu}$ is extended by zero outside of $D$ in the
following theorems.

\medskip

\begin{corol}\label{12corDIR1} {\it Let $D$ be a bounded domain in
${\Bbb C}$ whose boundary consists of $n\geqslant2$ mutually
disjoint Jordan curves and $\mu:D\to{\Bbb C}$ be a measurable
function  such that
\begin{equation}\label{eqDIR6*}\overline{\lim\limits_{\varepsilon\to0}}\quad
\dashint_{B(z_0,\varepsilon)}K_{\mu}(z)\,dm(z)<\infty\qquad\forall\
z_0\in\overline{D}\,.\end{equation} Then the Beltrami equation
(\ref{eqBeltrami}) has a multi-valued solutions of the Dirichlet
problem (\ref{eqGrUsl}) for each continuous function
$\varphi:\partial D\to{\Bbb R}$.} \end{corol}

\medskip

\begin{theo}\label{12thDIR2io} {\it Let $D$ be a bounded domain in
${\Bbb C}$ whose boundary consists of $n\geqslant2$ mutually
disjoint Jordan curves and $\mu:D\to\Bbb{C}$ be a measurable
function with $|\mu(z)|<1$ a.e. If $K_{\mu}\in L^1_{\rm loc}(D)$
and satisfies the condition
\begin{equation}\label{eq8.11.2}\int\limits_{0}^{\delta(z_0)}
\frac{dr}{||K_{\mu}||_{1}(z_0,r)}=\infty\end{equation} for some
$\delta(z_0)\in(0,d(z_0))$ where $d(z_0)=\sup\limits_{z\in
D}|z-z_0|$ and \begin{equation}\label{eq8.11.4}
||K_{\mu}||_{1}(z_0,r)=\int\limits_{D\cap S(z_0,r)}K_{\mu}(z)\,
|dz|\,,\end{equation} at each point $z_0\in\overline{D}$. Then the
Beltrami equation (\ref{eqBeltrami}) has a multi-valued solutions of
the Dirichlet problem (\ref{eqGrUsl}) for each continuous function
$\varphi:\partial D\to{\Bbb R}$.}\end{theo}

\medskip

\begin{corol}\label{12corDIR2t} {\it Let $D$ be a bounded domain in
${\Bbb C}$ whose boundary consists of $n\geqslant2$ mutually
disjoint Jordan curves and $\mu:D\to{\Bbb C}$ be a measurable
function such that \begin{equation}\label{eqDIR61}
k_{z_0}(\varepsilon)=O\left(\log\frac{1}{\varepsilon}\right)\quad
\text{as}\quad\varepsilon\to0\quad\forall\
z_0\in\overline{D}\end{equation} where $k_{z_0}(\varepsilon)$ is
the average of the function $K_{\mu}(z)$ over
$S(z_0,\varepsilon)$. Then the Beltrami equation
(\ref{eqBeltrami}) has a multi-valued solutions of the Dirichlet
problem (\ref{eqGrUsl}) for each continuous function
$\varphi:\partial D\to{\Bbb R}$.}\end{corol}

\medskip

\begin{rem}\label{rem2} In particular, the conclusion holds if
\begin{equation}\label{eqDIR6**}K_{\mu}(z)=O\left(\log\frac{1}{|z-z_0|}\right)\qquad{\rm
as}\quad z\to z_0\quad\forall\
z_0\in\overline{D}\,.\end{equation}\end{rem}

\medskip

\begin{theo}\label{12thKR4.1} {\it Let $D$ be a bounded domain in
${\Bbb C}$ whose boundary consists of $n\geqslant2$ mutually
disjoint Jordan curves and $\mu:D\to{\Bbb C}$ be a measurable
function such that $|\mu(z)|<1$ a.e. and
\begin{equation}\label{eqKR4.1}\int\limits_{D}\Phi\left(K_{\mu}(z)\right)\,dm(z)<\infty\end{equation}
for a convex non-decreasing function
$\Phi:[0,\infty]\to[0,\infty]$. If \begin{equation}\label{eqKR4.2}
\int\limits_{\delta}^{\infty}\frac{d\tau}{\tau\Phi^{-1}(\tau)}=\infty\end{equation}
for some $\delta>\Phi(0)$. Then the Beltrami equation
(\ref{eqBeltrami}) has a multi-valued solutions of the Dirichlet
problem (\ref{eqGrUsl}) for each continuous function
$\varphi:\partial D\to{\Bbb R}$.} \end{theo}

\medskip

\begin{corol}\label{12corDIR1000} {\it In particular, the conclusion of
Theorem \ref{12thKR4.1} holds if, for some $\alpha>0$,
\begin{equation}\label{eqKR4.1}\int\limits_{D}e^{\alpha K_{\mu}(z)}\,dm(z)<\infty\,.\end{equation}}
\end{corol}

\medskip

\noindent Denis Kovtonyuk, Igor Petkov, Vladimir Ryazanov, Ruslan Salimov, \\
Institute of Applied Mathematics and Mechanics,\\
National Academy of Sciences of Ukraine, \\
74 Roze Luxemburg str., 83114 Donetsk, UKRAINE \\
denis$\underline{\ }$\,kovtonyuk@bk.ru, igorpetkov@i.ua, \\
vlryazanov1@rambler.ru, vl.ryazanov1@gmail.com, salimov07@rambler.ru
\end{document}